\documentclass[12pt]{amsart}
\usepackage{amsthm,amsmath,amstext,amsfonts,amssymb,color}

\newcommand{\Q}{{\mathbb Q}}
\newcommand{\Z}{{\mathbb Z}}
\newcommand{\R}{{\mathbb R}}
\renewcommand{\P}{{\mathbb P}}
\newcommand{\F}{{\mathcal F}}
\newcommand{\T}{{\tau}}

\newcommand{\meas}{{\mathbb \mu}}

\newtheorem{thm}{Theorem}
\newtheorem{lemma}[thm]{Lemma}

\newcommand{\dof}{\bf}

\begin{document}
\title{A Nonmeasurable Set from Coin Flips}

\author{Alexander E. Holroyd}
\author{Terry Soo}

\address{Department of Mathematics, University of British Columbia,
121 -- 1984 Mathematics Rd, Vancouver, BC V6T 1Z2, Canada.}
\email{holroyd@math.ubc.ca; tsoo@math.ubc.ca}

\maketitle
To motivate the elaborate machinery of measure theory, it is desirable to show that in some natural space $\Omega$ one cannot define a measure on {\em all} subsets of $\Omega$, if the measure is to satisfy certain natural properties.   The usual example is given by the Vitali set, obtained by
choosing one representative from each equivalence class of $\R$
induced by the relation $x \sim y$ if and only if $x - y \in \Q$.  The resulting set is not meas\-ur\-able with respect to any translation-invariant measure on $\R$ that gives non-zero, finite measure to the unit interval \cite{vitali}.  In particular, the resulting set is not Lebesgue meas\-ur\-able.  The  construction above uses the axiom of choice.  Indeed, the Solovay theorem
\cite{Solovay} states that in the absence of the axiom of choice, there is a model of Zermelo-Frankel set theory where all the subsets of $\R$ are Lebesgue meas\-ur\-able.

In this note we give a variant proof of the existence of a
non\-meas\-ur\-able set (in a slightly different space).
We will use the axiom of choice in the guise of
the well-ordering principle (see the later discussion for more
information).  Other examples of non\-meas\-ur\-able sets may be found for example
in \cite{blackwell} and \cite[Ch.\ 5]{oxtoby}.

We will produce a non\-meas\-ur\-able set in the space $\Omega:=\{0,1\}^{\Z}.$  Translation-invariance plays a key role in the Vitali proof.  Here shift-invariance will play a similar role.    The
{\dof shift} $T: \Z \to \Z$ on  integers is defined via $Tx:=x+1$, and the shift $\T:\Omega \to \Omega$ on elements
$\omega\in\Omega$ is defined via $(\T\omega)(x):=\omega(x-1).$  We write $\T A:= \{\T\omega : \omega \in A\}$ for $A \subseteq \Omega.$

\begin{thm}
\label{main}
Let $\mathcal F$ be a $\sigma$-algebra on $\Omega$ that contains all singletons and is closed under the shift (that is, $A \in \mathcal F$ implies $\T A \in \mathcal F$).  If there exists a measure $\mu$ on $\mathcal F$ that is shift-invariant (that is, $\mu = \mu \circ \T$) and satisfies $\mu(\Omega) \in (0,\infty)$, and $\mu(\{\omega\}) = 0$ for all $\omega \in \Omega$, then $\mathcal F$ does not contain all subsets of $\Omega$.
\end{thm}

The conditions on $\mathcal{F}$ and $\mu$ in Theorem {\ref{main}} are indeed satisfied by measures that arise naturally.  A central example is the probability space $(\Omega, \mathcal G, \P)$ for a sequence of independent fair coin flips indexed by $\Z$, which is defined as follows.
Let $\mathcal A$ be the algebra of all
sets of the form $\{ \omega \in \Omega : \omega(k) = a_k, \text{for all} \ k \in K\}$, where $K \subset \Z$ is any finite subset of the integers
and $a \in \{0,1\}^K$ is any finite binary string.  The measure $\P$
restricted to $\mathcal A$ is given by $\P \big( {\{ \omega \in \Omega
:\omega(k) = a_k, \ \text{for all} \ k \in K\}} \big) = 2^{-|K|}$, where $|K|$
denotes the cardinality of $K$.  Thus $\P(\Omega) =1$, and $\P = \P\circ \T$ on $\mathcal A$.  The  Carath\'{e}odory extension theorem {\cite[Ch.\ 12, Theorem 8]{royden} gives a unique extension $\P$ to $\mathcal{G}:=\sigma(\mathcal A)$ (the $\sigma$-algebra generated by A) satisfying $\P = \P \circ \T$.  In addition, the continuity of measure implies $\P(\{\omega\}) = 0$ for all $\omega \in \Omega$.  Hence Theorem \ref{main} implies that $\mathcal G$ does not contain all subsets of $\Omega$.  Of course, the same holds for any extension $(\Omega, \mathcal G', \P')$ of $(\Omega, \mathcal G, \P)$ for which $\P'$ is shift-invariant (such as the completion under $\P$).

To prove Theorem {\ref{main}} we will define a non\-meas\-ur\-able function.  We are interested in functions from $\Omega$ to $\Z$ that are defined
everywhere except on some set of measure zero.  Therefore, for
convenience, introduce an additional element $\Delta \not \in \Z$.
Consider a function $X:\Omega\to\Z\cup\{\Delta\}.$ We call $X$ {\dof
almost-everywhere defined} if $X^{-1}\{\Delta\}$ is countable, which implies that $\meas(X^{-1}\{\Delta\}) = 0,$ for any measure $\mu$ satisfying the conditions of Theorem {\ref{main}}.  A function $X$ is {\dof meas\-ur\-able} with respect to $\F$ if $X^{-1}\{x\} \in \F$ for all $x \in \Z$.  We call $X$ {\dof shift-equi\-vari\-ant} if
$$X(\T\omega)=T(X(\omega)) \quad\text{for all }\omega\in\Omega$$
(where $T(\Delta):=\Delta$).  (We may think of a shift-equivariant $X$ as an ``origin-independent'' rule for choosing an element from the sequence $\omega$.)  Shift-equivariant functions of random processes are important in many settings, including percolation theory (for example in \cite{burtonkeane}) and coding theory (for example in \cite{keanea,keaneb}).

\begin{lemma}
\label{not-exist}
If $X:\Omega\to\Z\cup\{\Delta\}$ is an almost-everywhere defined, shift-equi\-vari\-ant function then $X$ is not meas\-ur\-able with respect to any $\F$ satisfying the conditions of Theorem {\ref{main}}.
\end{lemma}
\begin{lemma}
\label{exist}
There exists an almost-everywhere defined, shift-equi\-vari\-ant function
$X:\Omega\to\Z\cup\{\Delta\}$.
\end{lemma}

Theorem \ref{main} is an immediate consequence of the preceding two facts.
\begin{proof}[Proof of Theorem {\ref{main}}]
Let $(\Omega, \F, \meas)$ be a measure space satisfying the conditions of Theorem {\ref{main}}.  Using Lemma {\ref{exist}}, let $X$ be an almost-every\-where defined
shift-equi\-vari\-ant function.  By Lemma {\ref{not-exist}}, $X$ is not
$\F$-meas\-ur\-able.  Therefore there exists $z \in \Z$
such that $X^{-1}\{z\} \not \in \F.$
\end{proof}

\begin{proof}[Proof of Lemma \ref{not-exist}]
Towards a contradiction, let $X$ be a meas\-ur\-able function on $(\Omega, \F, \meas)$ satisfying the
conditions of Lemma {\ref{not-exist}}.
Since $X$ is shift-equi\-vari\-ant we have for each $x\in\Z$,
$$\meas(X^{-1}\{x\})=\meas\big(\T^{-x}X^{-1}\{x\}\big)=\meas(X^{-1}\{0\}).$$
Hence
$$\meas(X^{-1}\Z) = \meas\Big(\bigcup_{x\in\Z} X^{-1}\{x\}\Big) =\sum_{x\in
\Z} \meas(X^{-1}\{0\}) =0\text{ or }\infty,$$ which contradicts the facts that
$\meas(X^{-1}\{\Delta\}) = 0$ and $\meas(\Omega) \in (0, \infty)$.
\end{proof}

Let us recall some
facts about well-ordering.  A total order $\preceq$ on a set $W$
is a {\dof{well order}} if every nonempty subset of $W$ has a least
element.  The well-ordering principle states that every set has a
well order.  It is a classical result of Zermelo
{\cite{well-ordering}} that the well-ordering principle is equivalent
to the axiom of choice.
\begin{proof}[Proof of Lemma \ref{exist}]
Say $\omega \in \Omega$ is {\dof periodic} if $\T^x\omega=\omega$ for some $x\in\Z
\setminus \{0\}$.  If $\omega$ is not periodic then
$(\T^x\omega)_{x\in\Z}$ are all distinct.    Using the well-ordering
principle, fix a well order $\preceq$ of $\Omega$ and define the
function
$$X(\omega):=
\begin{cases}
\Delta &\text{if $\omega$ is periodic;}\\ {\displaystyle \text{the
unique} \ x \ \text{minimizing} \ \T^{-x}\omega \ \text{under} \
\preceq} & \text{otherwise.}\\
\end{cases}
$$
(We may think of $\T^{-x}\omega$ as $\omega$ viewed from location $x$, in which case $X$ is the location from which $\omega$ appears least.)
Clearly, $X$ is shift-equi\-vari\-ant.  It is
almost-everywhere defined since $\Omega$ contains only countably many periodic
elements.
\end{proof}

\paragraph{Acknowledgments}Alexander E. Holroyd is funded in part by an NSERC (Canada) Discovery Grant.  Terry Soo is funded in part by an NSERC PGS D and a UBC Graduate fellowship.

\providecommand{\bysame}{\leavevmode\hbox to3em{\hrulefill}\thinspace}
\providecommand{\MR}{\relax\ifhmode\unskip\space\fi MR }
\providecommand{\MRhref}[2]{%
  \href{http://www.ams.org/mathscinet-getitem?mr=#1}{#2}
}
\providecommand{\href}[2]{#2}


\begin{thebibliography}{1}

\bibitem{blackwell}
D.~Blackwell and P.~Diaconis, A non-measurable tail set, in \emph{Statistics,
  {P}robability and {G}ame {T}heory}, IMS Lecture Notes-Monograph Series,
  vol.~30, Institute of Mathematical Statistics, Hayward, CA, 1996, 1--5.

\bibitem{burtonkeane}
R.~M. Burton and M.~Keane, Density and uniqueness in percolation, \emph{Comm.
  Math. Phys.} \textbf{121} (1989) 501--505.

\bibitem{keanea}
M.~Keane and M.~Smorodinsky, A class of finitary codes, \emph{Israel J. Math.}
  \textbf{26} (1977) 352--371.

\bibitem{keaneb}
\bysame, Bernoulli schemes of the same entropy are finitarily isomorphic,
  \emph{Ann. of Math. (2)} \textbf{109} (1979) 397--406.

\bibitem{oxtoby}
J.~C. Oxtoby, \emph{Measure and Category}, 2nd ed., Graduate Texts in
  Mathematics, vol.~2, Springer-Verlag, New York, 1980.

\bibitem{royden}
H.~L. Royden, \emph{Real Analysis}, 3rd ed., Macmillan, New York, 1988.

\bibitem{Solovay}
R.~M. Solovay, A model of set-theory in which every set of reals is {L}ebesgue
  measurable, \emph{Ann. of Math. (2)} \textbf{92} (1970) 1--56.

\bibitem{vitali}
G.~Vitali, Sul problema della misura dei gruppi di punti di una retta,
  Gamberini and Parmeggiani, Bologna, 1905.

\bibitem{well-ordering}
E.~Zermelo, Beweis, da\ss\ jede {M}enge wohlgeordnet werden kann, \emph{Math.
  Ann.} \textbf{59} (1904) 514--516.

\end{thebibliography}
\end{document}